\documentclass[11pt, twoside]{article}
\usepackage{latexsym}
\usepackage{amsmath}
\usepackage{amssymb}
\usepackage[all]{xy}
\usepackage{amsfonts}
\usepackage{verbatim}
\usepackage{amsthm}
\usepackage{mathrsfs}
\usepackage{epsfig}
\usepackage{xy}
\usepackage{array}
\usepackage{stmaryrd}
\usepackage{graphicx,color}
\usepackage{xcolor}
\usepackage{tikz}
\usetikzlibrary{arrows,calc}
\usepackage{etex}
\usepackage{mathdots}
\usepackage{float}
\usepackage{graphics}
\usepackage{pdflscape}
\usepackage{xcolor}
\usepackage[colorlinks=true,linkcolor=blue,citecolor=blue]{hyperref}
\usepackage{anysize,hyperref}
\input xypic
\xyoption{all}

\usepackage[perpage,symbol]{footmisc}
\usepackage{setspace}
\def\E{\mathbb{E}}
\def\s{\mathfrak{s}}

\def\del{\delta}
\def\dr{\ar@{->}[r]}

\newcommand{\tri}[7]{\xymatrix@C=1.5em{#1\ar[r]^{#5}&#2\ar[r]^{#6}&#3\ar[r]^{#7}&#4}}
\newcommand{\ltri}[8]{\xymatrix@C=#8cm{#1\ar[r]^{#5}&#2\ar[r]^{#6}&#3\ar[r]^{#7}&#4}}
\newcommand{\rtri}[6]{\xymatrix@C=1.5em{#1\ar[r]^{#4}&#2\ar[r]^{#5}&#3\ar[r]^{#6}&\Sigma #1}}

\begin{document}
\baselineskip=15pt
\title{\Large{\bf Schanuel's lemma for extriangulated categories}}
\medskip
\author{Hangyu Yin}

\date{}

\maketitle
\def\blue{\color{blue}}
\def\red{\color{red}}

\newtheorem{theorem}{Theorem}[section]
\newtheorem{lemma}[theorem]{Lemma}
\newtheorem{corollary}[theorem]{Corollary}
\newtheorem{proposition}[theorem]{Proposition}
\newtheorem{conjecture}{Conjecture}
\theoremstyle{definition}
\newtheorem{definition}[theorem]{Definition}
\newtheorem{question}[theorem]{Question}
\newtheorem{remark}[theorem]{Remark}
\newtheorem{remark*}[]{Remark}
\newtheorem{example}[theorem]{Example}
\newtheorem{example*}[]{Example}
\newtheorem{condition}[theorem]{Condition}
\newtheorem{condition*}[]{Condition}
\newtheorem{construction}[theorem]{Construction}
\newtheorem{construction*}[]{Construction}

\newtheorem{assumption}[theorem]{Assumption}
\newtheorem{assumption*}[]{Assumption}

\baselineskip=17pt
\parindent=0.5cm
\vspace{-6mm}

\begin{abstract}
\baselineskip=16pt
In the context of extriangulated categories, we establish the injective version of Schanuel's lemma in homological algebra.\\[4mm]
\textbf{Keywords}: extriangulated category; homological algebra; injective dimension\\[1mm]
\textbf{2020 Mathematics Subject Classification:} 18G80; 18E10
\end{abstract}

\pagestyle{myheadings}
\markboth{\rightline {\scriptsize   Hangyu Yin }}
         {\leftline{\scriptsize Schanuel's lemma for extriangulated categories}}

\section{Introduction}

In the realm of mathematics, specifically within the field of algebra, particularly module theory, Schanuel's lemma, attributed to Stephen Schanuel, offers a means to assess the extent to which modules deviate from being projective. This lemma proves valuable in establishing the Heller operator within the stable category and providing straightforward explanations for dimension shifting phenomena, see \cite{Ka,La}.

In \cite{Bu}, B\"{u}hler developed the field of homological algebra within the framework of Quillen's exact categories, with a particular emphasis on bounded cohomology.
Mathieu and Rosbotham \cite{MR} elucidated the workings of Schanuel's lemma in the context of general exact categories and established the Injective Dimension Theorem in full detail.

Extriangulated categories have been introduced relatively recently by Nakaoka and Palu \cite{NP}. They have been defined by isolating specific properties of ${\rm Ext}^1(-,-)$ within exact categories and triangulated categories, focusing on their relevance in the context of cotorsion pairs. Notably, both triangulated categories and exact categories can be categorized as extriangulated categories. However, it's worth mentioning that there exist numerous examples of extriangulated categories that do not fall under the category of either triangulated or exact categories, see \cite{NP,ZZ, NP1, FHZZ}.

In this note, we offer a comprehensive explanation of the mechanics of Schanuel's lemma in extriangulated categories. Furthermore, we present the rigorous establishment of the Injective Dimension Theorem, see Theorem \ref{main1}.

\section{Preliminaries}

Let us briefly recapitulate some definitions and fundamental properties of extriangulated categories as presented in \cite{NP}. While we omit certain details here, readers can find them in \cite{NP}.

Consider an additive category $\mathcal{C}$ equipped with an additive bifunctor
\[
\mathbb{E}: \mathcal{C}^{\mathrm{op}} \times \mathcal{C} \rightarrow \mathrm{Ab},
\]
where $\mathrm{Ab}$ denotes the category of abelian groups. For any objects $A, C \in \mathcal{C}$, an element $\delta \in \mathbb{E}(C,A)$ is called an $\mathbb{E}$-extension. Let $\mathfrak{s}$ be a correspondence that associates an equivalence class
\[
\mathfrak{s}(\delta) = \xymatrix@C=0.8cm{[A \ar[r]^x & B \ar[r]^y & C]}
\]
with any $\mathbb{E}$-extension $\delta \in \mathbb{E}(C, A)$. This $\mathfrak{s}$ is referred to as a \textit{realization} of $\mathbb{E}$ if it renders the diagrams in \cite[Definition 2.9]{NP} commutative. A triplet $(\mathcal{C}, \mathbb{E}, \mathfrak{s})$ is called an \textit{extriangulated category} when it satisfies the following conditions:
\begin{enumerate}
\item[(1)] $\mathbb{E}\colon\mathcal{C}^{\mathrm{op}} \times \mathcal{C} \rightarrow \mathrm{Ab}$ is an additive bifunctor.

\item[(2)] $\mathfrak{s}$ is an additive realization of $\mathbb{E}$.

\item[(3)] $\mathbb{E}$ and $\mathfrak{s}$ satisfy some axioms (ET3), (ET3)$^{\rm op}$,
(ET4) and (ET4)$^{\rm op}$ in \cite[Definition 2.12]{NP}.
\end{enumerate}

We compile the following terminology as outlined in \cite{NP}.

\begin{definition}
Let $(\mathcal{C},\E,\s)$ be an extriangulated category.
\begin{enumerate}
\item[(1)] A sequence $A\xrightarrow{~x~}B\xrightarrow{~y~}C$ is referred to as a \textit{conflation} if it realizes an $\mathbb{E}$-extension $\delta\in\mathbb{E}(C,A)$. In this scenario, $x$ is called an \textit{inflation}, and $y$ is called a \textit{deflation}. We introduce two classes of morphisms: $\mathscr{M}$ consisting of all inflations, and $\mathscr{P}$ comprising all deflations.

\item[(2)] If a conflation $A\xrightarrow{~x~}B\xrightarrow{~y~}C$ realizes $\delta\in\mathbb{E}(C,A)$, we denote the pair $(A\xrightarrow{~x~}B\xrightarrow{~y~}C,\delta)$ as an \textit{$\mathbb{E}$-triangle}, presented as follows:
\[
A\overset{x}{\longrightarrow}B\overset{y}{\longrightarrow}C\overset{\delta}{\dashrightarrow}
\]
We usually omit the "$\delta$" if it is not required in the context.

\item[(3)] Let $A\overset{x}{\longrightarrow}B\overset{y}{\longrightarrow}C\overset{\delta}{\dashrightarrow}$ and $A^{\prime}\overset{x^{\prime}}{\longrightarrow}B^{\prime}\overset{y^{\prime}}{\longrightarrow}C^{\prime}\overset{\delta^{\prime}}{\dashrightarrow}$ be any pair of $\mathbb{E}$-triangles. If a triplet $(a,b,c)$ realizes $(a,c)\colon\delta\to\delta^{\prime}$, we represent it as:
\[
\xymatrix{
A \ar[r]^x \ar[d]^a & B\ar[r]^y \ar[d]^{b} & C\ar@{-->}[r]^{\delta}\ar[d]^c&\\
A'\ar[r]^{x'} & B' \ar[r]^{y'} & C'\ar@{-->}[r]^{\delta'} &}
\]
This triplet $(a,b,c)$ is termed a \textit{morphism of $\mathbb{E}$-triangles}.

\item[(4)] An object $I\in\mathcal{C}$ is considered \textit{injective} if, for any $\mathbb{E}$-triangle $\xymatrix{E\ar[r]^{\mu}&F\ar[r]^{\pi}&G\ar@{-->}[r]^{\delta}&}$ and any morphism $f\in\mathcal{C}(E,I)$, there exists $g\in\mathcal{C}(F,I)$ satisfying $g\mu=f$. The full subcategory containing injective objects in $\mathcal{C}$ is denoted as $\mathrm{inj}(\mathcal{C})$. Similarly, the full subcategory of projective objects is denoted as $\mathrm{proj}(\mathcal{C})$.

\item[(5)] We say that $\mathcal{C}$ \textit{has enough injectives} if, for any object $E\in\mathcal{C}$, there exists an $\mathbb{E}$-triangle:
\[
\xymatrix{E\ar[r]^{x}&I\ar[r]^{y}&F\ar@{-->}[r]^{\delta}&}
\]
where $I$ belongs to $\mathrm{inj}(\mathcal{C})$. A dual definition for having enough projectives can be formulated.
\end{enumerate}
\end{definition}

We provide some equivalent characterizations of injective objects, which will be used later.

\begin{proposition}\label{prop22}
Let $(\mathcal{C},\mathbb{E},\s)$ be an extriangulated category and $E$ be an object in $\mathcal{C}$.
The following statements are equivalent:

{\rm (i)} $E$ is injective;

{\rm (ii)} Every $\E$-triangle $\xymatrix{E\ar[r]^{\mu}&F\ar[r]^{\pi}&G\ar@{-->}[r]^{\delta}&}$
in $\mathcal{C}$, where $\mu$ is a split monomorphism;

{\rm (iii)} There exists an $\E$-triangle $\xymatrix{E\ar[r]^{\mu}&I\ar[r]^{\pi}&G\ar@{-->}[r]^{\delta}&}$ satisfying $I$ is injective and $\mu$ is a split monomorphism.
\end{proposition}

\proof  (i) $\Rightarrow$ (ii) For any $\E$-triangle $\xymatrix{E\ar[r]^{\mu}&F\ar[r]^{\pi}&G\ar@{-->}[r]^{\delta}&}$ in $\mathcal{C}$, by (i), we know that
$E$ is injective, thus there exists $\mu :F\rightarrow E$, such that $\tilde{\mu}\mu=id_E$.
This shows that $u$ is a split monomorphism.

(ii) $\Rightarrow$ (i) For any $\E$-triangle $\xymatrix{A\ar[r]^{a}&F\ar[r]^{b}&B\ar@{-->}[r]^{\delta}&}$ and $f:A \rightarrow E$,
by \cite[Proposition 1.20]{LN}, there exists a morphism of $\E$-triangles
$$\xymatrix{
A \ar[r]^a \ar[d]^f & F\ar[r]^b \ar[d]^{g} & B\ar@{-->}[r]^{\del}\ar@{=}[d]&\\
E\ar[r]^{a'} & F' \ar[r]^{b'} & B\ar@{-->}[r]^{\del'} &.}$$
By (ii), $a'$ has a left inverse $c$, then $f=cga$.
This shows that $E$ is injective.

(i) $\Rightarrow$ (iii) Let $E$ is injective. For $\E$-triangle $\xymatrix{E\ar[r]^{id_E}&E\ar[r]&0\ar@{-->}[r]^{0}&}$,
$\xymatrix{E\ar[r]^{id_E}&E}$ is a split monomorphism by (ii).

(iii) $\Rightarrow$ (ii) For any $\E$-triangle $\xymatrix{E\ar[r]^{\mu}&F\ar[r]^{\pi}&G\ar@{-->}[r]^{\delta}&,}$ there exist an $\E$-triangle $$\xymatrix{E\ar[r]^{i}&I\ar[r]^{p}&G'\ar@{-->}[r]^{\delta'}&}$$
 satisfying $I$ is injective and $\mu$ is a split monomorphism by (iii). Since $I$ is  injective, then there exist $\alpha\colon F\rightarrow I$ such that $\alpha \mu=i$. Since $i$ has a left inverse $i'$, then $\mu$ has a left inverse $i'\alpha$.
 This shows that $\mu$ is a split monomorphism.
 \qed
\vspace{2mm}

The following proposition can be found in \cite[Proposition 2.3]{MR}.

\begin{proposition}\label{prop23}
 Let  $E, F, G$  be three objects in an additive category $\mathscr{C}$. The following
 statements are equivalent:

{\rm (i)} $ F $ is a product of $ E$  and $ G $;

{\rm (ii)}  $F $ is a coproduct of  $E$  and $ G $;

{\rm (iii)} There exists an $\E$-triangle $$\xymatrix{E\ar[r]^{\mu}&F\ar[r]^{\pi}&G\ar@{-->}[r]^{\delta}&}$$
and morphisms  $\tilde{\mu} \in  \mathcal{C}(F, E)$  and  $\tilde{\pi} \in \mathcal{C}(G, F)$  such that  $\tilde{\mu} \circ \mu=\operatorname{id}_{E}$ and $\pi \circ \tilde{\pi}=\mathrm{id}_{G}$, and $\mu \circ \tilde{\mu}+\tilde{\pi} \circ \pi=\mathrm{id}_{F}$;

{\rm (iv)} There exists an $\E$-triangle $$\xymatrix{E\ar[r]^{\mu}&F\ar[r]^{\pi}&G\ar@{-->}[r]^{\delta}&}$$
and a morphism $ \tilde{\mu} \in \mathcal{C}(F, E)$  such that  $\tilde{\mu} \circ \mu=\mathrm{id}_{E}$, where $\mathrm{id}_{E}$ the identity morphism on  $E$.
\end{proposition}

\begin{proposition}\label{prop24}
Assume that
$\xymatrix{E\ar[r]^{\mu}&F\ar[r]^{\pi}&G\ar@{-->}[r]^{\delta}&}$ is an $\E$-triangle.
If  $F\cong E \oplus G$, then  $F$  is injective if and only if both  $E$  and $ G$  are injective.
\end{proposition}

\proof
For any $\E$-triangle $\xymatrix{A\ar[r]^{a}&B\ar[r]^{b}&C\ar@{-->}[r]^{\delta}&}$ and $f:A\rightarrow E$, there exists $\alpha \in \mathcal{C}(B,F)$ such that $\alpha a=\mu f$ since $f$   is injective. Because $F \cong E \oplus G$, then there exists $\tilde{\mu}\in\mathcal{C}(F,E)$ such that $\tilde{\mu}\mu={\rm id}_E$. It follows that $\tilde{\mu}\alpha a=f$. Therefore $E$ is injective. Similarly,  we show that $G$ is injective. \qed

\section{Schanuel's lemma for extriangulated categories}

To establish the main results of this section, we require the following preparatory work.
Unless otherwise specified, we assume that $(\mathcal{C},\mathbb{E},\s)$ is an extriangulated category from now on.

\begin{lemma}
Let  $(\mathcal{C}, \mathbb{E}, \mathfrak{s})$  be an extriangulated category. Let
$$\xymatrix{A\ar[r]^{h}&C\ar[r]^{h'}&E\ar@{-->}[r]^{\delta_h}&,}~
\xymatrix{D\ar[r]^{d}&E\ar[r]^{d'}&F\ar@{-->}[r]^{\delta_d}&,}~
\xymatrix{B\ar[r]^{g}&C\ar[r]^{g'}&F\ar@{-->}[r]^{\delta_g}&}$$
be three $\mathbb{E}$-triangles satisfying $ h=g \circ f$. Then there exists commutative diagram
$$\xymatrix{
A \ar[r]^f \ar@{=}[d] & B\ar[r]^{f'} \ar[d]^{g} & D\ar[d]^d &\\
A\ar[r]^{h} & C \ar[r]^{h'} \ar[d]^{g'} & E \ar[d]^{d'} & \\
& F \ar@{=}[r] & F &}$$
which satisfies

{\rm (i)} $\xymatrix{A\ar[r]^{f}&B\ar[r]^{f'}&D\ar@{-->}[r]^{d^*\delta_h}&}$ is an $\mathbb{E}$-triangle;

{\rm (ii)} $f'_*\delta_g=\delta_d$;

{\rm (iii)} $f_{*}\delta_{h}=d'^{*}\delta_{g}$.
\end{lemma}

\proof By (ET4)$^{\rm  op}$, there exist an $\E$-triangle $\xymatrix{B_0\ar[r]^{g_0}&C\ar[r]^{g'}&F\ar@{-->}[r]^{\delta_{g_0}}&}$
and a commutative diagram in $\mathcal{C}$
$$\xymatrix{
A \ar[r]^{f_0} \ar@{=}[d] & B_0\ar[r]^{f_0'} \ar[d]^{g_0} & D\ar[d]^d &\\
A\ar[r]^{h} & C \ar[r]^{h'} \ar[d]^{g'} & E \ar[d]^{d'} & \\
& F \ar@{=}[r] & F &}$$
which satisfies

$\rm (i')$  $ \xymatrix{A\ar[r]^{f_0}&B_0\ar[r]^{f_0'}&D\ar@{-->}[r]^{d^*\delta_h}&}$ is an  $\mathbb{E}$-triangle,

$\rm (ii')$  ${f_0'}_*\delta_{g_0}=\delta_d,$

$\rm (iii')$  ${f_0}_{*}\delta_{h}=d'^{*}\delta_{g_0} $.

By (ET3)$^{\rm op}$, we obtain a morphism of $\E$-triangles
$$\xymatrix{
B_0\ar[r]^{g_0} \ar@{-->}[d]^u & C\ar[r]^{g'} \ar@{=}[d] & F\ar@{-->}[r]^{\del_{g_0}}\ar@{=}[d]&\\
B\ar[r]^{g} & C \ar[r]^{g'} & F \ar@{-->}[r]^{\del_g} &.}$$
We obtain that $u$ is an isomorphism by \cite[Corollary 3.6]{NP}.
In particular, we have $u_*\del_{g_0}=\del_g$. If we put $f=uf_0$ and $f'=f_0'u^{-1}$, then we obtain a commutative diagram
$$\xymatrix{
A \ar[r]^f \ar@{=}[d] & B\ar[r]^{f'} \ar[d]^{g} & D\ar[d]^d &\\
A\ar[r]^{h} & C \ar[r]^{h'} \ar[d]^{g'} & E \ar[d]^{d'} & \\
& F \ar@{=}[r] & F &}$$

$\rm (i)$ By \cite[Proposition 3.7]{NP}, we have $[\xymatrix{A\ar[r]^f & B\ar[r]^{f'} & D}]=[\xymatrix{A\ar[r]^{f_0} & B_0\ar[r]^{f_0'} & D}]$, then $ \xymatrix{A\ar[r]^{f}&B\ar[r]^{f'}&D\ar@{-->}[r]^{d^*\delta_h}&} $ is an  $\mathbb{E} $-triangle;

$\rm (ii)$ By ${\rm (ii')}$ ${f_0'}_*\del_{g_0}=\del_d$, we have ${f'}_*\del_g=(f_0'u^{-1})_*\del_g={f_0'}_*u^{-1}_*\del_g
={f_0'}_*\del_{g_0}=\del_d;$

$\rm (iii)$ By ${\rm (iii')}$ ${f_0}_{*}\delta_{h}=d'^{*}\delta_{g_0} $,
we have
$$f_{*}\delta_{h}=(uf_0)_*\del_h=u_*{f_0}_*\del_h
=u_*{d'}^*\delta_{g_0}={d'}^*u_*\delta_{g_0}={d'}^*\del_g.$$ \qed

\begin{lemma}\label{lem32}
Let  $(\mathcal{C}, \mathbb{E}, \mathfrak{s})$ be an extriangulated category.
Let $\xymatrix{A\ar[r]^{x_1}&B_1\ar[r]^{y_1}&C_1\ar@{-->}[r]^{\delta_1}&}$
and
$\xymatrix{A\ar[r]^{x_2}&B_2\ar[r]^{y_2}&C_2\ar@{-->}[r]^{\delta_2}&}$
be two $\mathbb{E}$-triangles. Then there exists a commutative diagram
$$\xymatrix{A\ar[r]^{x_1} \ar[d]^{x_2} & B_1\ar[r]^{y_1} \ar[d]^{m_2} & C_1\ar@{=}[d]\\
B_2\ar[r]^{m_1} \ar[d]^{y_2} & M\ar[r]^{e_1} \ar[d]^{e_2}&C_1\\
C_2\ar@{=}[r] & C_2 &}
$$
which satisfies:

${\rm (i)}~\mathfrak{s}(x_{2}^{*} \delta_{1})=[\xymatrix{B_2\ar[r]^{m_1}&M\ar[r]^{e_1}&C_1}];$

${\rm (ii)}~\mathfrak{s}(x_{1 *} \delta_{2})=[\xymatrix{B_1\ar[r]^{m_2}&M\ar[r]^{e_2}&C_2}];$

${\rm (iii)}~ e_{1}^* \delta_{1}+e_{2}^* \delta_{2}=0.$
\end{lemma}
\proof
By the additivity of  $\mathfrak{s},$ we have  $\mathfrak{s}\left(\delta_{1} \oplus \delta_{2}\right)=[\xymatrix{A \oplus A\ar[r]^{x_{1} \oplus x_{2}~}&B_{1} \oplus B_{2}\ar[r]^{y_{1} \oplus y_{2}}&C_1 \oplus C_2}].$
Let  $C_{1} \underset{p_{1}}{\stackrel{\iota_{1}}{\rightleftarrows}} C_{1} \oplus C_{2} \underset{p_{2}}{\stackrel{\iota_{2}}{\leftrightarrows}} C_{2}$  be a biproduct in  $\mathcal{C}$. Put  $\mu=\left(\triangledown_{C}\right)_{*}\left(\delta_{1} \oplus \delta_{2}\right)$  and take $$\mathfrak{s}(\mu)=
[\xymatrix{A\ar[r]^{j}&M \ar[r]^{k~~~}&C_1 \oplus C_2}].$$ Then  $\mu $ satisfies
$$\iota_1^* \mu=\delta_{1} \text { and } \iota_2^* \mu=\delta_{2} \text {. }$$
Applying (ET4)$^{\rm op}$ to  $\mathfrak{s}(0)=[\xymatrix@C=0.6cm{C_1\ar[r]^{\iota_{1}~~~}&C_1\oplus C_{2} \ar[r]^{~~~p_2}&C_2}]$ and $\mathfrak{s}(\mu)=
[\xymatrix@C=0.6cm{A\ar[r]^{j}&M \ar[r]^{k~~~}&C_1 \oplus C_2}],$  there exists an $\E$-triangle $\xymatrix{B_1'\ar[r]^{m_2'}&M\ar[r]^{e_2}&C_2\ar@{-->}[r]^{\theta_1}&} $  and a commutative diagram
$$\xymatrix{
A \ar[r]^{x_1'} \ar@{=}[d] & B_1'\ar[r]^{y_1'} \ar[d]^{m_2'} & C_1\ar[d]^{\iota_1} &\\
A\ar[r]^{j} & M \ar[r]^{k} \ar[d]^{e_2} & C_1\oplus C_2 \ar[d]^{p_2} & \\
& C_2 \ar@{=}[r] & C_2 &}$$
which satisfies

${\rm (i')}$ $ \xymatrix{A\ar[r]^{x_1'}&B_1'\ar[r]^{y_1'}&C_1\ar@{-->}[r]^{\iota_1^*\mu}&} $ is an  $\mathbb{E} $-triangle;

${\rm (ii')}$  $y'_{1*}\theta_1=0$;

${\rm (iii')}$  $x'_{1*}\mu=p_2^*\theta_1$.

Moreover $[\xymatrix{A\ar[r]^{x_1'} & B_1'\ar[r]^{y_1'}&C_1}]=\s(\iota_1^*\mu)=[\xymatrix{A\ar[r]^{x_1} & B_1 \ar[r]^{y_1} &C_1}].$

Thus there is an isomorphism  $b_{1} \in \mathcal{C}(B_{1}, B_{1}')$ satisfying $b_{1} \circ x_{1}=x_{1}'$ and  $y_{1}' \circ b_{1}=y_{1}.$ If we put $ m_{2}=m_{2}' \circ b_1$, then we obtain that $ \xymatrix{B_1\ar[r]^{m_2}&M\ar[r]^{e_2}&C_2\ar@{-->}[r]^{b_{1*}^{-1}\theta_1}&} $ is an $\E$-triangle by \cite[Proposition 3.7]{NP}.
Put $\theta=b_{1*}^{-1}\theta_1$,
then we obtain a commutative diagram
$$\xymatrix{
A \ar[r]^{x_1} \ar@{=}[d] & B_1\ar[r]^{y_1} \ar[d]^{m_2} & C_1\ar[d]^{\iota_1} &\\
A\ar[r]^{j} & M \ar[r]^{k~~~} \ar[d]^{e_2} & C_1\oplus C_2 \ar[d]^{p_2} & \\
& C_2 \ar@{=}[r] & C_2 &}$$
which satisfies $$x_{1*}\delta_2=x_{1*}\iota_2^*\mu=(b_1^{-1}x_1')_*\iota_2^*\mu
=b_{1*}^{-1}\iota_2^*x_{1*}'\mu=b_{1*}^{-1}\iota_2^*p_2^*\theta_1
=b_{1*}^{-1}(p_2\iota_2)^*\theta_1=b_{1*}^{-1}\theta_1=\theta.$$
Thus we obtain $ \mathfrak{s}(x_{1*}\del_2)=\mathfrak{s}(\theta)
=[\xymatrix{B_1\ar[r]^{m_2}&M\ar[r]^{e_2}&C_2}]$.

Similarly, from $\s(0)=[\xymatrix{C_2\ar[r]^{\iota_2~~~}&C_1\oplus C_2\ar[r]^{~~~p_2}&C_1}]$ and $\s(\mu)=[\xymatrix{A\ar[r]^j&M\ar[r]^{k~~~}&C_1\oplus C_2}]$,
we obtain a commutative diagram
$$\xymatrix{
A \ar[r]^{x_2} \ar@{=}[d] & B_2\ar[r]^{y_2} \ar[d]^{m_1} & C_2\ar[d]^{\iota_2} &\\
A\ar[r]^{j} & M \ar[r]^{k~~~} \ar[d]^{e_1} & C_1\oplus C_2 \ar[d]^{p_1} & \\
& C_1 \ar@{=}[r] & C_1 &}$$
which satisfies $\s(x_{2*}\del_1)=[\xymatrix{B_2\ar[r]^{m_1}&M\ar[r]^{e_1}&C_1}]$. Since $$m_1x_2=j=m_2x_1,~e_1m_2=p_1km_2=p_1\iota_1y_1=y_1,~e_2m_1=p_2km_1=p_2\iota_2y_2=y_2.$$
Moreover, we have $e_1^*\del_1+e_2^*\del_2=k^*p_1^*\del_1+k^*p_2^*\del_2
=k^*(p_1^*\del_1+p_2^*\del_2)=k^*(p_1^*\iota_1^*\mu+p_2^*\iota_2^*\mu)
=k^*((\iota_1 p_1)^*\mu+(\iota_2 p_2)^*\mu)=k^*((\iota_1 p_1+\iota_2 p_2)^*\mu)=k^*\mu=0$ by Lemma \ref{lem32}. \qed

\begin{proposition}\label{prop33}

Suppose that $\xymatrix{E\ar[r]^{\mu}&I\ar[r]^{\pi}&F\ar@{-->}[r]^{\delta}&}$
and $\xymatrix{E\ar[r]^{\mu'}&I'\ar[r]^{\pi'}&F\ar@{-->}[r]^{\delta'}&}$
are two $\E$-triangles, and $I$, $I'$ are injective. Then $I \oplus F'\cong I' \oplus F$.
\end{proposition}

\proof
For two $\E$-triangles $\xymatrix{E\ar[r]^{\mu}&I\ar[r]^{\pi}&F\ar@{-->}[r]^{\delta}&}$
and $\xymatrix{E\ar[r]^{\mu'}&I'\ar[r]^{\pi'}&F\ar@{-->}[r]^{\delta'}&}$, we obtain that a commutative diagram by Lemma \ref{lem32}
$$\xymatrix{E\ar[r]^{\mu} \ar[d]^{\mu'} & I\ar[r]^{\pi} \ar[d]^{h} & F\ar@{=}[d]\\
I'\ar[r]^{h'} \ar[d]^{\pi'} & C\ar[r]^{p} \ar[d]^{p'}&F\\
F'\ar@{=}[r] & F' &}
$$
which satisfies $\xymatrix{I\ar[r]^h& C\ar[r]^{p'}&F'\ar@{-->}[r]^{\mu_*\del'}&}$
and $\xymatrix{I'\ar[r]^{h'}& C\ar[r]^{p}&F\ar@{-->}[r]^{\mu_*'\del}&}$
are two $\E$-triangles.
Since $I$, $I'$ are injective, then $\mu$ and $\mu'$ are split monomorphisms.
 Therefore $C \cong I\oplus F'$ and $C\cong I'\oplus F$ by Proposition \ref{prop23}.
  It follows that $I\oplus F'\cong I'\oplus F$. \qed
\vspace{2mm}

This proposition immediately yields the following conclusion.

\begin{corollary}\label{cor34}
Let $(\mathcal{C}, \E ,\s)$ an extriangulated category.
Suppose that there exists a diagram of morphisms in $(\mathcal{C}, \E ,\s)$ of the form
$$\xymatrix{E\ar[r]\ar[d]^{\cong}&I\ar[r]&F\ar@{-->}[r]&\\
E'\ar[r]&I'\ar[r]&F'\ar@{-->}[r]&}$$
such that $I$ and $I'$ are injective, the horizontal lines are $\E$-triangles and the vertical arrow is an isomorphism. Then $I\oplus F'\cong I'\oplus F$.

\end{corollary}

\begin{definition}
For an object $E \in \mathcal{C}$, an injective resolution of $E$ is a sequence of admissible morphisms of the form:
\[
\xymatrix@C=0.5cm{
E&&I^0&&\cdots &&I^{n-1}&&I^n&&\cdots\\
&G^0~&&G^1~&&G^{n-1}~~&&G^n~&&G^{n+1}~~\\
\ar"1,1";"1,3"
\ar"1,3";"1,5"
\ar"1,5";"1,7"
\ar"1,7";"1,9"
\ar"1,9";"1,11"
\ar@{->>}_{\cong}"1,1";"2,2"
\ar@{->>}"1,3";"2,4"
\ar@{->>}"1,5";"2,6"
\ar@{->>}"1,7";"2,8"
\ar@{->>}"1,9";"2,10"
\ar@{>->}"2,2";"1,3"
\ar@{>->}"2,4";"1,5"
\ar@{>->}"2,6";"1,7"
\ar@{>->}"2,8";"1,9"
\ar@{>->}"2,10";"1,11"
}
\]
such that, for each $n \geq 0$, the object $I^{n}$ is injective, and
$\xymatrix{G^n~\ar@{>->}[r]&I^n\ar@{->>}[r]&G^{n+1}}$
is a conflation. It is obvious that if $\mathcal{C}$ has enough injectives, then we can build an injective resolution for every object in $\mathcal{C}$.
\end{definition}

\begin{proposition}
Suppose that we have the following injective resolutions of $E$, with the factorization of each admissible morphism included:
\[
\xymatrix@C=0.5cm{
E&&I^0&&\cdots &&I^{n-1}&&I^n&&\cdots\\
&G^0~&&G^1~&&G^{n-1}~~&&G^n~~&&G^{n+1}~~\\
\ar"1,1";"1,3"
\ar"1,3";"1,5"
\ar"1,5";"1,7"
\ar"1,7";"1,9"
\ar"1,9";"1,11"
\ar@{->>}_{\cong}"1,1";"2,2"
\ar@{->>}"1,3";"2,4"
\ar@{->>}"1,5";"2,6"
\ar@{->>}"1,7";"2,8"
\ar@{->>}"1,9";"2,10"
\ar@{>->}"2,2";"1,3"
\ar@{>->}"2,4";"1,5"
\ar@{>->}"2,6";"1,7"
\ar@{>->}"2,8";"1,9"
\ar@{>->}"2,10";"1,11"
}
\]
and
\[
\xymatrix@C=0.5cm{
E&&J^0&&\cdots &&J^{n-1}&&J^n&&\cdots\\
&H^0~&&H^1~&&H^{n-1}~~&&H^n~~&&H^{n+1}~~\\
\ar"1,1";"1,3"
\ar"1,3";"1,5"
\ar"1,5";"1,7"
\ar"1,7";"1,9"
\ar"1,9";"1,11"
\ar@{->>}_{\cong}"1,1";"2,2"
\ar@{->>}"1,3";"2,4"
\ar@{->>}"1,5";"2,6"
\ar@{->>}"1,7";"2,8"
\ar@{->>}"1,9";"2,10"
\ar@{>->}"2,2";"1,3"
\ar@{>->}"2,4";"1,5"
\ar@{>->}"2,6";"1,7"
\ar@{>->}"2,8";"1,9"
\ar@{>->}"2,10";"1,11"
}
\]
Then, for each $n \geq 1$, we have isomorphisms
$$I^{0} \oplus J^{1} \oplus I^{2} \oplus \cdots \oplus J^{2 n-1} \oplus G^{2 n} \cong J^{0} \oplus I^{1} \oplus J^{2} \oplus \cdots \oplus I^{2 n-1} \oplus H^{2 n}$$
and
$$I^{0} \oplus J^{1} \oplus I^{2} \oplus \cdots \oplus J^{2 n-1} \oplus I^{2 n} \oplus H^{2 n+1} \cong J^{0} \oplus I^{1} \oplus J^{2} \oplus \cdots \oplus I^{2 n-1} \oplus J^{2 n} \oplus G^{2 n+1}.$$

\proof We prove this by induction. For $n=1$, first note that Corollary \ref{cor34}, applied to the diagram
$$\xymatrix{G^0\ar[r]\ar[d]^{\cong}&I_0\ar[r]&G^1\ar@{-->}[r]&\\
H^0\ar[r]&J^0\ar[r]&H^1\ar@{-->}[r]&}$$
obtains  $I^{0} \oplus H^{1} \cong J^{0} \oplus G^{1}$. Further, there is a diagram of the form
$$\xymatrix{I^0\oplus H^1\ar[r]\ar[d]^{\cong}&I_0\oplus J^1\ar[r]&H^2\ar@{-->}[r]&\\
J^0\oplus G^1\ar[r]&J^0\oplus I^1\ar[r]&G^2\ar@{-->}[r]&.}$$
By Corollary \ref{cor34}, we have  $I^{0} \oplus J^{1} \oplus G^{2} \cong J^{0} \oplus I^{1} \oplus H^{2}$. To finish the proof for $n=1$, we again get a diagram
$$\xymatrix{I^{0} \oplus J^{1} \oplus G^{2}\ar[r]\ar[d]^{\cong}&I^{0} \oplus J^{1} \oplus I^{2}\ar[r]&G^{3}\ar@{-->}[r]&\\
J^{0} \oplus I^{1} \oplus H^{2}\ar[r]&J^{0} \oplus I^{1} \oplus J^{2}\ar[r]&H^{3}\ar@{-->}[r]&.}$$
It follows that, by Corollary \ref{cor34} we obtain an isomorphism $ I^{0} \oplus J^{1} \oplus I^{2} \oplus H^{3} \cong J^{0} \oplus I^{1} \oplus J^{2} \oplus G^{3}.$

Assume the result holds some $ n \geq 1$. Then there exists a diagram

$$\xymatrix{I^{0} \oplus \cdots \oplus I^{2 n} \oplus H^{2 n+1}\ar[r]\ar[d]^{\cong}
&I^{0} \oplus \cdots \oplus I^{2 n} \oplus J^{2 n+1}\ar[r]
&H^{2(n+1)}\ar@{-->}[r]&\\
J^{0} \oplus \cdots \oplus J^{2 n} \oplus G^{2 n+1}\ar[r]
&J^{0} \oplus \cdots \oplus J^{2 n} \oplus I^{2 n+1}\ar[r]
&G^{2(n+1)}\ar@{-->}[r]&.}$$
By Corollary \ref{cor34}, we have
$$\begin{array}{l}
I^{0} \oplus J^{1} \oplus I^{2} \oplus \cdots \oplus J^{2(n+1)-1} \oplus G^{2(n+1)} \\
\quad \cong J^{0} \oplus I^{1} \oplus J^{2} \oplus \cdots \oplus I^{2(n+1)-1} \oplus H^{2(n+1)} .
\end{array}$$
There exists a diagram:
$$\xymatrix@C=0.5cm{I^{0} \oplus J^{1} \oplus I^{2} \oplus \cdots \oplus J^{2n+1} \oplus G^{2(n+1)}\ar[r]\ar[d]^{\cong}&I^{0} \oplus J^{1} \oplus I^{2} \oplus \cdots \oplus J^{2n+1} \oplus I^{2(n+1)}\ar[r]&G^{2(n+1)+1}\ar@{-->}[r]&\\
J^{0} \oplus I^{1} \oplus J^{2} \oplus \cdots \oplus I^{2n+1} \oplus H^{2(n+1)}\ar[r]&J^{0} \oplus I^{1} \oplus J^{2} \oplus \cdots \oplus I^{2n+1} \oplus J^{2(n+1)}\ar[r]&H^{2(n+1)+1}\ar@{-->}[r]&.}$$
By Corollary \ref{cor34}, we have
$$\begin{array}{l}
I^{0} \oplus J^{1} \oplus I^{2} \oplus \cdots \oplus J^{2(n+1)-1} \oplus I^{2(n+1)} \oplus H^{2(n+1)+1} \\
\quad \cong J^{0} \oplus I^{1} \oplus J^{2} \oplus \cdots \oplus I^{2(n+1)-1} \oplus J^{2(n+1)} \oplus G^{2(n+1)+1}.
\end{array}$$
The conclusion is proven through induction. \qed

\end{proposition}

\begin{theorem}\label{main1}{\bf (Injective Dimension Theorem)}
Let $(\mathcal{C},\E,\s)$ be an extriangulated category. Suppose that $\mathscr{C}$  has enough injectives. The following statements are equivalent for  $n \geq 1$ and every  $E \in \mathcal{C}$:

{\rm (i)} If there is an exact sequence of admissible morphisms
$$\xymatrix{E~~\ar@{>->}[r]&I^0\ar[r]
&\cdots\ar[r]&I^{n-1}\ar@{->>}[r]&F} \eqno (1)$$
with each $I^{m}$, $0 \leq m \leq n-1$ injective, then $F$ must be injective;

{\rm (ii)} There is an exact sequence of admissible morphisms
$$\xymatrix{E~~\ar@{>->}[r]&I^0\ar[r]
&\cdots\ar[r]&I^{n-1}\ar@{->>}[r]&I^n} \eqno (2)$$
with each  $I^{m}$, $0 \leq m \leq n$ injective.

\end{theorem}

\proof ${\rm (i)} \Rightarrow {\rm (ii)}.$ Since  $\mathcal{C}$  has enough injectives, we can build an injective resolution of $E$:
\[
\xymatrix@C=0.6cm{
E&&I^0&&\cdots &&I^{n-1}&&J^n&\cdots\\
&G^0~&&G^1~&&G^{n-1}~~&&I^n~\\
\ar"1,1";"1,3"
\ar"1,3";"1,5"
\ar"1,5";"1,7"
\ar"1,7";"1,9"
\ar"1,9";"1,10"
\ar@{->>}_{\cong}"1,1";"2,2"
\ar@{->>}"1,3";"2,4"
\ar@{->>}"1,5";"2,6"
\ar@{->>}"1,7";"2,8"
\ar@{>->}"2,2";"1,3"
\ar@{>->}"2,4";"1,5"
\ar@{>->}"2,6";"1,7"
\ar@{>->}"2,8";"1,9"
}
\eqno (3)\]
We obtain an exact sequence as in diagram (2), and $I^{n}$ must be injective, by the condition {\rm (i)}.

${\rm (ii)} \Rightarrow {\rm (i)}$. There must are an injective resolution of  $E$ of the form
\[
\xymatrix@C=0.6cm{
E&&J^0&&\cdots &&J^{n-1}&&J^n&\cdots\\
&J^0~&&J^1~&&J^{n-1}~~&&J^n~\\
\ar"1,1";"1,3"
\ar"1,3";"1,5"
\ar"1,5";"1,7"
\ar"1,7";"1,9"
\ar"1,9";"1,10"
\ar@{->>}_{\cong}"1,1";"2,2"
\ar@{->>}"1,3";"2,4"
\ar@{->>}"1,5";"2,6"
\ar@{->>}"1,7";"2,8"
\ar@{>->}"2,2";"1,3"
\ar@{>->}"2,4";"1,5"
\ar@{>->}"2,6";"1,7"
\ar@{>->}_{\cong}"2,8";"1,9"
}
\]
and for any exact sequence as in diagram (1), with each  $I^{n}$  injective, there exists an injective resolution
\[
\xymatrix@C=0.6cm{
E&&I^0&&\cdots &&I^{n-1}&&I^n&\cdots\\
&G^0~&&G^1~&&G^{n-1}~~&&G^n~\\
\ar"1,1";"1,3"
\ar"1,3";"1,5"
\ar"1,5";"1,7"
\ar"1,7";"1,9"
\ar"1,9";"1,10"
\ar@{->>}_{\cong}"1,1";"2,2"
\ar@{->>}"1,3";"2,4"
\ar@{->>}"1,5";"2,6"
\ar@{->>}"1,7";"2,8"
\ar@{>->}"2,2";"1,3"
\ar@{>->}"2,4";"1,5"
\ar@{>->}"2,6";"1,7"
\ar@{>->}"2,8";"1,9"
}
\]
with $G^{n}=F$.

When $n$ is an odd number, put $$I=I^{0} \oplus J^{1} \oplus I^{2} \oplus \cdots \oplus I^{n-1} \oplus J^n,~~G=J^{0} \oplus I^{1} \oplus J^{2} \oplus \cdots \oplus J^{n-1},$$ then $I\cong G\oplus F$;

When $n$ is an even number, put $$G=I^{0} \oplus J^{1} \oplus I^{2} \oplus \cdots \oplus J^{n-1},~~I=J^{0} \oplus I^{1} \oplus J^{2} \oplus \cdots \oplus I^{n-1} \oplus J^n,$$ then $I\cong G\oplus F$.
Therefore there exists an $\E$-triangle $$\xymatrix{F\ar[r]^{\mu}&I\ar[r]^{\pi}&G\ar@{-->}[r]^{\delta}&}$$
such that $\mu$ is a split monomorphism. Since $I$ is a finite product of injective objects, then $I$ is injective by Proposition \ref{prop24}. It follows that by Proposition \ref{prop22} we obtain that $F$ is injective.   \qed

\begin{definition}
Let $\mathscr{M}$ be the class of inflations in an extriangulated category $(\mathcal{C},\E,\s)$. We say that $E\in \mathcal{C}$ has finite injective dimension if there exists an exact sequence of admissible morphisms as shown in diagram (2), where each $I^{m}$ is an injective object. If $E$ has finite injective dimension, we write $\operatorname{Inj}{\mathscr{M}}\text{-dim}(E)=0$ if $E$ is injective, and $\operatorname{Inj}{\mathscr{M}}\text{-dim}(E)=n$ if $E$ is not injective and $n$ is the smallest natural number for which there exists an exact sequence of admissible morphisms as shown in diagram (2) where every $I^{m}$ is injective. If $E$ does not have finite injective dimension, we write $\operatorname{Inj}_{\mathscr{M}}\text{-dim}(E)=\infty$.

The global dimension of the extriangulated category $(\mathcal{C},\E,\s)$ is
$$\sup \left\{\operatorname{Inj}_{\mathscr{M}}\text{-dim}(E) \mid E \in \mathcal{C}\right\} \in \mathbb{N}_{0} \cup\{\infty\}.$$
\end{definition}

\begin{remark}
The injective dimension of an object $E$ in an extriangulated category $(\mathcal{C},\E,\s)$ can be determined by examining any of its injective resolutions. Indeed, consider the following as an injective resolution of $E$ (with the factorization of each admissible morphism included):
\[
\xymatrix@C=0.6cm{
E&&J^0&&\cdots &&J^{n-1}&&J^n&\cdots\\
&G^0~&&G^1~&&G^{n-1}~~&&G^n~\\
\ar"1,1";"1,3"
\ar"1,3";"1,5"
\ar"1,5";"1,7"
\ar"1,7";"1,9"
\ar"1,9";"1,10"
\ar@{->>}_{\cong}"1,1";"2,2"
\ar@{->>}"1,3";"2,4"
\ar@{->>}"1,5";"2,6"
\ar@{->>}"1,7";"2,8"
\ar@{>->}"2,2";"1,3"
\ar@{>->}"2,4";"1,5"
\ar@{>->}"2,6";"1,7"
\ar@{>->}"2,8";"1,9"
}
\]
Then, by Theorem \ref{main1}, ${\rm Inj}_{\mathscr{M}}\text{-dim}(E)\leq n $ if and only if $G^{n}$ is injective.

\end{remark}

\textbf{Hangyu Yin}\\
College of Mathematics, Hunan Institute of Science and Technology, 414006, Yueyang, Hunan, P. R. China.\\
E-mail: hangyuyin@163.com

\end{document}